\documentclass[11pt]{article}
\usepackage{amsfonts,amsmath,amstext,amsbsy,euscript,amssymb, graphicx}

\makeatletter
 \def\LaTeX{\leavevmode L\raise.42ex
   \hbox{\kern-.3em\size{\sf@size}{0pt}\selectfont A}\kern-.15em\TeX}
\makeatother

\newcommand{\BibTeX}{{\rm B\kern-.05em{\sci\kern-.025emb}\kern-.08em\TeX}}

\newtheorem{thm}{Theorem}[section]
\newtheorem{lem}[thm]{Lemma}

\newtheorem{cor}[thm]{Corollary}

\newcommand{\us}{\underline{s}}

\newcommand{\Q}{\mathcal{Q}}

\numberwithin{equation}{section}

\begin{document}

\title{Lattice and Schr\"oder paths with periodic boundaries}

\author{Joseph P.S. Kung$^{1}$, Anna de Mier$^{2,5}$, Xinyu Sun$^3$,
and Catherine Yan$^{4,6}$  \\
\vspace{.3cm} \\
$^{1}$ Department of Mathematics, \\
University of North Texas,  Denton, TX 76203, U.S.A. and 
\\  Centre de Recerca Matem\`atica, 08193 Bellaterra, Spain \\
\vspace{.3cm} \\
$^{2}$ Department of Applied Mathematics 2, \\
Universitat Polit\`ecnica de Catalunya, Jordi Girona 1--3, 08034 Barcelona, Spain \\
\vspace{.3cm} \\
$^{3}$ Mathematics Department, \\ 
Tulane University, 6823 St. Charles Avenue, New Orleans, LA 70118,
U.S.A. \\
\vspace{.3cm} \\
$^{4}$ Department of Mathematics \\
Texas A\&M University,  College Station, TX 77843, U.S.A. and \\
Center for Combinatorics, LPMC-TJKLC, \\
Nankai University, Tianjin 300071,  P.R. China \\
 \vspace{.3cm} \\
$^1$ kung@unt.edu, $^2$ anna.de.mier@upc.edu,
$^3$xsun1@tulane.edu, $^4$cyan@math.tamu.edu}
\maketitle

\footnotetext[5]{The second author was supported by the ``Ram\'on y
  Cajal'' programme of the Spanish Ministry of Science and Technology.}
\footnotetext[6]{The fourth author was  supported in part by NSF grants 
DMS-0245526 and DMS-0653846.}


\begin{abstract}
We consider paths in the plane with $(1,0),$ $(0,1),$ and $(a,b)$-steps
that start at the origin, end at height $n,$ and stay strictly to the left 
of a given non-decreasing right boundary.  We show that if the boundary
is periodic and has slope at most $b/a,$ then the
ordinary generating function for the number of such paths ending at
height $n$ is algebraic.  Our argument is in two parts.  We
use a simple combinatorial decomposition to obtain an Appell relation or
``umbral'' generating function, in which the power $z^n$ is replaced by
a power series of the form $z^n \phi_n(z),$ where $\phi_n(0) = 1.$
Then we convert (in an explicit way) the umbral generating function to
an ordinary generating function by solving a system of linear
equations and a polynomial equation.  This conversion implies
that the ordinary generating function is algebraic.

\end{abstract}

\section{Paths in the plane}

A {\sl lattice path} $\pi$ is a path in the plane with two
kinds of steps: a unit north step $N$ or a unit east step $E.$
Enumerating lattice paths and related combinatorial objects
is closely related to calculating
several probability distributions used in non-parametric statistics.
See, for example, \cite{Mohanty79, HN1981}.   If
$x$ is a positive integer, a lattice path from the origin $(0,0)$
to the point $(x-1,n)$ can be coded by a length-$n$ non-decreasing
sequence $(x_0,x_1, \ldots,x_{n-1}),$ where $0 \le x_i \le x-1$
and $x_i$ is the $x$-coordinate of the $(i+1)$-st north step.
For example, the path $EENENNEE$ is coded by
$(2,3,3)$ and the horizontal path ending at $(x-1,0)$ is coded by the
empty sequence.

Let $\us$ be a non-decreasing sequence with positive integer terms
$s_0,s_1,s_2, \ldots,$ thought of as a {\sl (right) boundary.} A
lattice path from $(0,0)$ to $(x-1,n)$ is an {\sl $\us$-lattice
path} if $ x_i < s_i $ for $0 \le i \le n-1.$
(Note that we require the path to stay {\it strictly to the left} of the
boundary.)
If $x \ge
s_{n-1},$ then the number of $\us$-lattice paths from $(0,0)$
to $(x-1,n)$ does not depend on $x.$  Let $\mathrm{LP}_{n}(\us)$
be this common number, which equals the number of $\us$-lattice
paths from $(0,0)$ to $(s_n-1,n).$ For instance,
the sequence $\underline{i+1}$ with terms $1,2,3, \ldots $
 gives a one-step staircase boundary.  A
 $45^{\circ}$-rotation gives a bijection between
$(\underline{i+1})$-lattice paths and Dyck paths. Thus,
$\underline{s}$-lattice paths can be thought of as Dyck paths on a
``bumpy'' $x$-axis.

A boundary $\underline{s}$ is {\sl periodic of height $k$ and width $l$}
if there is a finite non-decreasing
sequence $\us_0$ with $k$ terms $s_0, s_1, \ldots, s_{k-1}$ such that
$s_{k-1} = l$ and $ s_m = s_r + ql, $
where $q = \lfloor m/k \rfloor,$ and
$r = m - qk.$  Put another way, $\us$ is the concatenation
$ \us_0, \us_0 + l, \us_0 + 2l, \us_0 + 3l, \ldots \,\,. $
A boundary is {\sl ultimately periodic} if it can be written as the concatenation 
$\us^{\prime}, p+ \us ,$ where $\us^{\prime}$ is a finite initial sequence ending 
in $p$ and $\us$ is a periodic boundary.

In this paper, we present an elementary proof (avoiding the use of
cathalitic variables and hence the kernel method)
of a slightly more general form of a theorem of de Mier and Noy \cite{deMier}.

\begin{thm} \label{TBP}
Let $\us$ be an ultimately periodic right boundary. Then
the ordinary generating function
$$
\sum_{n = 0}^{\infty}  \mathrm{LP}_{n}(\us) \, z^n
$$
is algebraic.
\end{thm}

The proof begins with a combinatorial
decomposition which first appeared in the enumeration of parking functions
\cite{KonW, KungYan03}.   This decomposition yields an
Appell relation, an ``umbral'' generating function with
the power $z^n$ replaced by a power series $z^n \phi_n (z).$
When $\us$ is ultimately periodic,
the Appell
relation can be converted to an ordinary generating function
by solving
a polynomial equation and a system of linear equations.
The way this conversion is done implies that
the ordinary generating function is algebraic.
The method of obtaining generating functions from an Appell relation 
or functional equation is an old one and dates back, at least, to P\'olya.  
Special cases of results in this paper can be find in \cite{BSS, Gessel, Merlini, Tamm} 
and other papers.  

Our method extends to paths allowing a third kind of steps. Let
$a$ and $b$ be non-negative integers. An $(a,b)$-path $\pi$ is a
path starting at the origin with three kinds of steps: $N,$ $E,$
and {\sl diagonal} steps $D$ which go from a point $(u,v)$ to the
point $(u+a,v+b).$
For example, the path $EEDDNENNEE$ is a $(2,1)$-path starting at $(0,0)$
and ending at $(9,5).$  We allow $a$ and $b$ to be any non-negative
integer. If $a = b = 0,$ then $D$ is a null step or ``pause''.
If $a =0$ and $b=1$ or $a=1$ and $b=0,$ then we
 have two kinds
of $N$ or $E$ steps.

An $(a,b)$-path from $(0,0)$ to $(x-1,n)$ with $d$ diagonal steps
can be coded by a pair of sequences. Let $\pi$ be
such a path. Then there are $n - d(b-1)$ steps that are either $N$
or $D$. Replace every diagonal step $D$ by a $(0,b)$-step $C$ to
obtain the {\sl compacted} path. The compacted path goes from
$(0,0)$ to $(x-da-1,n).$   Record the $x$-coordinates where the
$N$ or $C$ steps occur.  This gives a non-decreasing sequence
$(x_0, x_1, \ldots, x_{n-d(b-1)-1})$ of length $n - d(b-1)$ with terms satisfying $0
\le x_i \le x - da -1.$  We also record the indices of the $C$
steps in this sequence.  This is an increasing sequence
$(y_0, y_1, \ldots, y_{d-1})$ of length $d$ with terms satisfying 
$0 \le y_i \le n-d(b-1)-1.$ For our example, the compacted path is $EECCNENNEE$
and $((2,2,2,3,3),(0,1))$ codes the original path $\pi.$

We may also think of an $(a,b)$-path $\pi$ as a {\sl polygon} or union of line
segments. Inductively, if the path ends at the point $(u,v),$ an
east step adds the unit horizontal line segment from $(u,v)$ to
$(u+1,v),$ a north step adds the unit vertical line segment from
$(u,v)$ to $(u,v+1),$ and a diagonal step adds the straight line segment
from $(u,v)$ to $(u+a,v+b)$ (with slope $b/a$).
The {\sl right extent}
$\mathrm{Right}_i(\pi)$ of
an $(a,b)$-path $\pi$ at $y$-coordinate $i$ is the $x$-coordinate of
the rightmost point of the intersection of the
path (as a polygon) with the horizontal line $y = i.$
An $(a,b)$-path from $(0,0)$ to $(x-1,n)$
is a {\sl $(a,b;\us)$-path} if for $0 \le i \le n-1,$
$$
\mathrm{Right}_i( \pi ) < s_i.
$$
As in the lattice path case, a $45^{\circ}$-rotation gives a
bijection between $(1,1;\underline{i+1})$-paths and Schr\"oder
paths.  

We begin this paper with a combinatorial decomposition for the set
of all paths in a rectangle (Section 2).  Using this decomposition,
we derive recursions and Appell relations in Section 3.
In Section 4, we show how to convert an Appell relation satisfying
a periodic condition into an ordinary generating function.  This method
gives a proof of a generalization
of Theorem \ref{TBP} for $(a,b)$-paths.  In Section 5, we work out
a concrete example, giving another solution to the tennis
ball problem.  We derive explicit formulas for $(1,b)$-paths with
 an arithmetic-progression boundary in Section 6.
 Section 7 is devoted to analogs of our results for parking
functions. We conclude, in Section 8, with some remarks and open problems.   

This paper is self-contained for the most part.  We shall only use
the elementary and formal part of the theory of algebraic fractional power series.
For background, see, for example, \cite[Chapter 6]{EC2}.

\section{Combinatorial decomposition}

In this section, $\mathrm{Path}_{n}(\us)$ denotes the set of all
$(a,b;\us)$-paths or the set of $\us$-lattice paths, depending on
the context, that end at $(s_n-1,n)$. If $y$ is a real number, let
$\underline{y}$ be the constant sequence with terms $y,y,y, \ldots.$

Given a non-negative integer $a$ and a positive integer $b,$
a boundary $\us$ satisfies the {\sl slope condition}
for $(a,b)$-paths if
for every $j,$
$$
s_{j+1} > s_j - 1 + a/b, \
s_{j+2} > s_j - 1 + 2a/b,  \ldots,\ 
s_{j+b} > s_j - 1 + a.
$$
We adopt the convention that the slope condition is satisfied for
lattice paths and $(0,0)$-paths, and not satisfied for $(a,0)$-paths if $a > 0.$   
Thinking of the boundary as the union of the line
segments from $(s_i,i)$ to $(s_{i+1},i+1),$ the slope condition
says that for every $j,$ the line segment $(s_j-1,j)$ to
$(s_j-1+a,j+b)$ (corresponding to a diagonal step) lies strictly
above the right boundary.  This implies that an $(a,b)$-path can
only reach a point $(s_i,i)$ by an east step. This technical
condition allows arguments for lattice paths to transfer without
change to $(a,b)$-paths but is reasonably natural.  (If the
slope condition is not assumed, the situation is more complicated.
However, there are analogs of our results for some cases.)

\begin{figure}[ht]
\begin{center}
\includegraphics{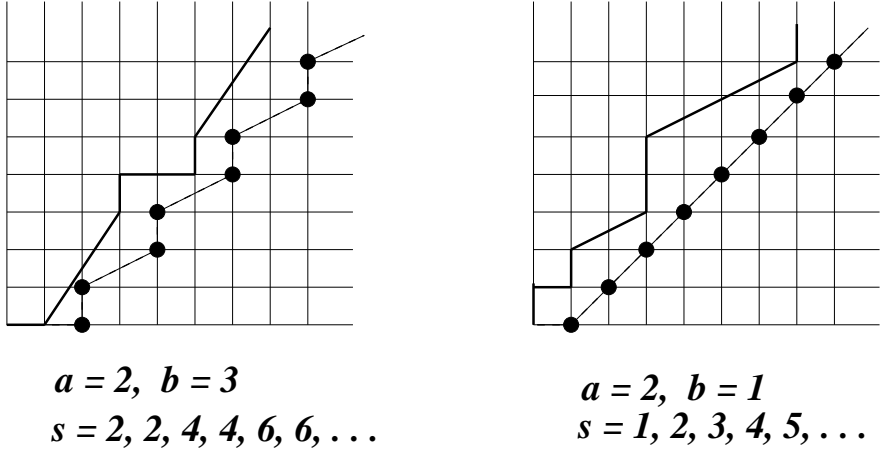}

\caption{Two examples of $(a,b;\us)$-paths. In the left-hand side
diagram, the border satisfies the slope condition, whereas the
border on the right-hand side does not.} \label{fig:us}

\end{center}

\end{figure}

For example, the boundary $\underline{i+1}$ satisfies the slope
condition for $(1,1)$-paths.  More generally, when $\gamma \le b/a,$
the boundary $\underline{\lceil i/\gamma \rceil + 1}$ of integers 
immediately to
the right of the line $y = \gamma (x-1),$ satisfies the slope
condition for $(a,b)$-paths.
See  Figure~\ref{fig:us} for two explicit examples.

\begin{thm} \label{decomposition}
Let $\us$ be a right boundary satisfying the slope
condition and $x$ be an integer such that
$x > s_{n}.$ Then
there is a bijection:
$$
\mathrm{Path}_n(\underline{x})
\longleftrightarrow
\bigcup_{m=0}^n
\mathrm{Path}_m ( \us ) \times
\mathrm{Path}_{n-m}(\underline{x-s_m}).
$$
\end{thm}

{\it Proof.} Observe that $\mathrm{Path}_n(\underline{x})$ is the
set of all paths in the rectangle with lower left corner $(0,0)$
and upper right corner $(x-1,n).$   For a path $\pi,$ let $m $ be
the (unique) index such that $\mathrm{Right}_i( \pi ) < s_i$ for
all $i < m$ and $\mathrm{Right}_m( \pi ) \ge  s_m.$ By the slope
condition, the path $\pi$ reaches the boundary by an east step.
Thus, every path in $\mathrm{Path}_n(\underline{x})$ can be
decomposed uniquely into two subpaths: a subpath in
$\mathrm{Path}_m (\us)$ and a subpath in
the rectangle with corners $(s_m,m)$ and $(x-1,n).$ These paths
are connected by an east step. To finish the proof, observe that
paths in the rectangle with corners $(s_m,m)$ and $(x-1,n)$ are in
bijection with paths in $\mathrm{Path}_{n-m}(\underline{x-s_m}).$
\hfill $\Box$

\section{Recursions and Appell relations}

The bijection of the previous section yields an infinite system of linear
recursions, which are then written as an Appell relation.  We first
consider lattice paths.  Since lattice paths are coded by
non-decreasing sequences or multisets,
$$
|\mathrm{Path}_n(\underline{x})| = \binom {x+n-1}{n} = (-1)^n \binom {-x}{n}.
$$
By Theorem \ref{decomposition}, we have
$$
\sum_{m=0}^n
(-1)^m \mathrm{LP}_m (\us)
\binom {s_m - x}{n-m}
=
\binom {-x}{n}.
$$
Since this recursion holds for all integers $x$ such that
$x > s_{n},$ it holds as a polynomial
identity in $x.$  Setting $x = 0,$ we obtain, $\mathrm{LP}_0 (\us) = 1,$
and for $n \ge 1,$
\begin{eqnarray} \label{LR_n}
\sum_{m=0}^n
(-1)^m  \mathrm{LP}_m (\us)
\binom {s_m}{n-m}
=
0.
\end{eqnarray} 

The next step is to multiply each of these equations by $t^n$ and sum over all
$n$. Recall the following (easy) lemma, which can be found, for
example, in \cite{KSY}. 

\begin{lem} \label{AbstractAppell}
Suppose that the sequence $g_n$ is the solution to the triangular system
\begin{eqnarray*}
b_{0,0} g_0 & = & a_0
\\
b_{1,1} g_1 + b_{1,0} g_0  & = & a_1
\\
b_{2,2} g_2 + b_{2,1} g_1 + b_{2,0} g_0  & = & a_2
\\
& \vdots &
\\
b_{m,m} g_m + b_{m,m-1} g_{m-1} + b_{m,m-2} g_{m-2} + \cdots + b_{m,1}g_1
+ b_{m,0} g_0 & = & a_m
\\
&\vdots&
\end{eqnarray*}
where for all $m,$ the diagonal coefficient $b_{m,m}$ is non-zero.
Then
\begin{eqnarray} \label{eqn-abs-Appell}
\sum_{m=0}^{\infty}  g_m t^m \phi_m(t)  =  \Psi(t),
\end{eqnarray}
where
\begin{eqnarray*}
\Psi (t)   =  \sum_{k=0}^{\infty}  a_k t^k \,\,\, \mathrm{and}
\,\,\, \phi_m(t)  =   \sum_{k=0}^{\infty} b_{m+k,m} t^k.
\end{eqnarray*}
\end{lem}


Applying Lemma \ref{AbstractAppell}, 
we obtain
$$
1
=
\sum_{n=0}^{\infty} \left[\sum_{k=n}^{\infty} \binom {s_n}{k-n} t^k  \right]
(-1)^n  \mathrm{LP}_n (\us)
=
\sum_{n=0}^{\infty}
(-1)^n \mathrm{LP}_n (\us)  t^n(1+t)^{s_n}.
$$
Changing variables from $t$ to $-t,$ we obtain an Appell relation for $\us$-lattice paths.

\begin{lem} \label{LatticeAppell}
For any boundary $\us$, the number of $\us$-lattice paths satisfy the equation
$$
\sum_{n=0}^\infty \mathrm{LP}_n(\us) t^n (1-t)^{s_n} = 1.
$$
\end{lem}

This equation for $\us$-lattice paths  has been obtained  by Gessel 
\cite{Gessel}. 
Appell relations for $(a,b)$-paths are more complicated.
We will enumerate $(a,b)$-paths
weighted by the number of diagonal steps.
Let
$$
\mathrm{SP}_n (a,b,\sigma;\us) =
\sum_{d=0}^{\lfloor n/b \rfloor}
 |\mathrm{Path}_{n,d}(\us)|  \sigma^d,
$$
where $\sigma$ is a variable and $\mathrm{Path}_{n,d}(\us)$ is the set
of
 $(a,b;\us)-$paths ending at
$(s_{n}-1,n)$ and having $d$ diagonal steps.  The (diagonal-step) enumerator is a polynomial in $\sigma$ 
except in the case $a = b = 0.$ 

From the coding given in the introduction,
the number of all $(a,b)$-paths with $d$ diagonal steps
from $(0,0)$ to $(x-1,n)$ is

$$
(-1)^{n - d(b-1)} \binom {-(x-da)}{n-d(b-1)} \binom {n-d(b-1)}{d}.
$$
Hence, the enumerator $\mathrm{SP}_n (a,b,\sigma;\underline{x})$
of such paths is given by
\begin{eqnarray*}
\mathrm{SP}_n (a,b,\sigma;\underline{x}) = \sum_{d=0}^{\lfloor n/b \rfloor}
(-1)^{n - d(b-1)}  \binom {da -x}{n - d(b-1)} \binom {n-d(b-1)}{d}
\sigma^d \, .
\end{eqnarray*}
Thus, by Theorem~\ref{decomposition}, for a right boundary $\us$ (satisfying
the slope condition for $(a,b)$-paths),
$$
\sum_{m=0}^n \mathrm{SP}_m (a,b,\sigma ; \us) \mathrm{SP}_{n-m}
(a,b,\sigma;\underline{x-s_m})
 = \mathrm{SP}_n (a,b,\sigma;\underline{x}).
$$
%
%
Setting $x = 0,$ we obtain the recursions \eqref{SRn}:
\begin{eqnarray} 
&& \sum_{m=0}^n \mathrm{SP}_m (a,b,\sigma ; \us)
\left(\sum_{j=0}^{\lfloor (n-m)/b \rfloor} (-1)^{n-m-j(b-1)}
\binom {s_m + ja}{n-m-j(b-1)} \binom {n-m-j(b-1)}{j} \sigma^j
\right) \nonumber
\\
&=&
\sum_{d = \lceil n/(a+b-1) \rceil}^{\lfloor n/b \rfloor}
(-1)^{n - d(b-1)}  \binom {da}{n - d(b-1)}
\binom {n-d(b-1)}{d} \sigma^d. 
\label{SRn}  
\end{eqnarray}
To get the smaller range for the right sum,
we use the fact that $\binom {c}{d} \neq 0$ if and only if $ c \ge d.$
We remark that determinantal formulas
(similar to formulas given in \cite{KungYan03})
for the enumerator of $(a,b)$-paths
can be obtained with Cramer's rule from the recursions \eqref{SRn}.

We use Lemma~\ref{AbstractAppell} again to express these recursions as an Appell relation.
We start with the classic case of Schr\"oder paths.
When $a = b = 1,$ the sum on the right side of ($\mathrm{SR}_n$) has one term,
$(-1)^n \sigma^n.$
Hence,
$$
\Psi(t) = \sum_{n=0}^{\infty} (-1)^n \sigma^n t^n = \frac {1}{1 + \sigma t}.
$$
On the left side, the quantity $\mathrm{SP}_m(1,1,\sigma; \us)$ occurs
multiplied by
\begin{eqnarray*}
(-1)^k \left(\binom {s_m}{k} \binom {k}{0} + \binom {s_m + 1}{k}
\binom {k}{1} \sigma
 + \binom {s_m + 2}{k} \binom {k}{2} \sigma^2  + \ldots
 + \binom {s_m + k}{k} \binom {k}{k} \sigma^k\right),
\end{eqnarray*}
where $k=n-m$.

Using the easy binomial coefficient identity
$$
\binom {s+i}{k} \binom {k}{i} =
\binom {s+i}{i} \binom {s}{k-i},
$$
we conclude that
\begin{eqnarray*}
\phi_m(t) &=& (1-t)^{s_m} \left[ 1 - \binom {s_m + 1}{1} \sigma
t + \binom {s_m + 2}{2} \sigma^2  t^2 - \binom {s_m + 3}{3}
\sigma^3 t^3 + \ldots \right]
\\
&=& (1-t)^{s_m} \left[ \sum_{i=0}^{\infty} \binom {-(s_m + 1)}{i}
\sigma^i  t^i \right]
\\
&=& (1-t)^{s_m} (1 + \sigma t)^{-(s_m + 1)}.
\end{eqnarray*}
This yields the Appell relation
$$
\sum_{n=0}^{\infty} \mathrm{SP}_n(1,1,\sigma; \us)  t^n \left( \frac
{1 - t}{1 + \sigma t} \right)^{s_n} = 1.
$$
Setting $\sigma=0,$ we recover the Appell relation for lattice paths
in Lemma \ref{LatticeAppell}.  Setting $\sigma=1,$ we obtain
an Appell relation for the total number of Schr\"oder paths.
Finally, setting $\sigma=-1,$ we obtain
$$
\sum_{n=0}^{\infty} \mathrm{SP}_n(1,1,-1; \us)  t^n = 1.
$$
This implies that when $n \ge 1,$ $\mathrm{SP}_n(1,1,-1; \us)  = 0,$
that is, the number of Schr\"oder paths from $(0,0)$
to $(n-1,n)$ with an even number
of diagonal steps equals the number of such Schr\"oder paths with an odd number
of diagonal steps for any boundary $\us$ satisfying the slope condition for $(1,1)$-paths.
For the case $s_i=i+1$,  this is known. See  \cite[Prop. 2.1]{BSS}.  The two
combinatorial proofs in \cite{BSS} proving the special case
can be applied without change to prove the general case.

By a similar argument, we obtain the following Appell relation
for $(1,b)$-paths:
$$
\sum_{n=0}^{\infty} \mathrm{SP}_n(1,b,\sigma; \us)  t^n \left( \frac
{1 - t}{1 + \sigma t^b} \right)^{s_n} = 1.
$$

For general $(a,b)$-paths,
we have
$$
\phi_m(t) = (1-t)^{s_m} \left[ \sum_{j=0}^{\infty} (-1)^j \binom {s_m +
aj}{j} (1-t)^{j(a-1)} \sigma^j  t^{bj} \right]
$$
and
$$
\Psi(t) = \sum_{d=0}^{\infty} (-1)^d \binom {da}{d}
(1-t)^{d(a-1)} (\sigma  t^b)^d.
$$
The series $\phi_m(t)$ and $\Psi(t)$  can be simplified by
the following lemma of Gould \cite{Gould}.

\begin{lem}\label{Gould}
For nonnegative integers $a$ and $s$,
$$
\sum_{j=0}^{\infty}
\binom {s+aj}{j} z^j= \frac{f^{s+1}}{a+(1-a)f},
$$
where
$$
f(z)= \sum_{m=0}^{\infty}  \frac {1}{am+1} \binom{am+1}{m}z^m,
$$
is the unique power series solution of the equation $ f(z) - 1=z[f(z)]^a. $
\end{lem}





There are at least two ways to prove Lemma \ref{Gould}.
Gould's proof in \cite{Gould}
begins by proving
$$
\sum_{n=0}^{\infty}  (-1)^n \binom {s + an}{n} t^n (1+t)^{s + (a-1)n}
= \frac {1}{1 + at},
$$
using the identities
$$
\sum_{i=0}^n (-1)^i \binom {n}{i} \binom {s + ai}{n} = (-1)^n a^n.  
$$
These identities can be proved by a finite difference argument.
The lemma then follows from the substitution $1 + t = 1/f$ (so
that $-t = (f-1)/f$), and  $z = (f-1)/f^a.$
Another way to prove Lemma~\ref{Gould} is to use the Lagrange inversion formula
to show that 
the coefficient of $z^j$ in
$$
\frac {f^{s+1}}{a+(1-a)f}
$$
equals
$$
\frac {1}{n} \sum_{i=0}^{j-1} (a-1)^{j-i-1} \binom{s + aj}{i} (ja-ia+s),
$$
and then show by induction on $j$ that this sum equals $\binom{s + aj}{j}.$

Using Lemma \ref{Gould}, we obtain
$$
\phi_n(t)=\frac{F(t)}{a+(1-a)F(t)} \left( (1-t)F(t)  \right)^{s_n}
$$
and
$$
\Psi(t)=\frac{F(t)}{a+(1-a)F(t)},
$$
where $ F(t) =f (-(1-t)^{a-1}t^b \sigma). $

Hence, for general $(a,b)$-paths we have the following Appell relation.

\begin{lem} \label{SPAppell}  Let $\underline{s}$ be a boundary
satisfying the slope condition for $(a,b)$-paths.  Then
$$
\sum_{n=0}^{\infty}  \mathrm{SP}_n (a,b,\sigma;\underline{s})
t^n [(1-t)F(t)]^{s_n} = 1,
$$
where $ F(t) =f (-(1-t)^{a-1}t^b \sigma) $ and $f(t)$ is the unique
power series solution of $f(t)-1 = tf(t)^a.$  
\end{lem}



In the next two sections, we show how to convert some Appell relations
to ordinary generating functions.  We motivate these results
with the simplest case, when the boundary is an arithmetic progression.
Let $s_i = c + id.$   For $(a,b)$-paths
with $d \ge a/b,$  the Appell relation can be written as
$$
\sum_{n=0}^{\infty}
\mathrm{SP}_n (a,b,\sigma;\underline{s}) z^n = \frac {1}{[(1-t)F(t)]^c},
$$
where $z = t[(1-t)F(t)]^d. $
Since $t$ is an algebraic function of $z$ (see Lemma~\ref{Solutions} below), we conclude that
the ordinary generating function
of $\mathrm{SP}_n (\sigma;\underline{s})$ is algebraic.
In the case of lattice paths, the
Appell relation simplifies to
\begin{eqnarray} \label{arith}
\sum_{n=0}^{\infty}
\mathrm{LP}_n (\underline{s}) z^n = \frac {1}{(1-t)^c},
\end{eqnarray}
where $ z = t(1-t)^d.$ From this, we obtain the explicit formula
$$
\mathrm{LP}_n (\underline{c+id}) =
\frac {c}{c + n(d+1)} \binom {c+n(d+1)}{n}.
$$
This formula is known when $c = 1$ (see, for example,
\cite[Section 1.5]{Mohanty79} or \cite[p.~175]{EC2}).  In particular, the series $f(z)$ in Lemma 3.3 is the generating function
$\sum_{m=0}^{\infty} \mathrm{LP}_{m} (\underline{1 + ai})z^m.$
One way to derive this formula is to apply the Lagrange inversion formula (see, for example, \cite[p.~42]{EC2}) 
to equation \eqref{arith}.
%
%
Another way is
to use a difference operator theory,
directly analogous to the differential operator theory for 
Gon\u{c}arov polynomials developed in \cite{KungYan03}.

\section{Periodicity implies algebraicity}

In this section, we show how to obtain an ordinary generating function from 
an Appell relation when a periodicity condition is
satisfied. We also show that this generating function is algebraic. We refer
the reader to~Section~6.1 of~\cite{EC2} for background on algebraic generating
functions. Recall that  a fractional
power series is a series of the form
$$
\sum_{i=0}^{\infty} c_i z^{i/N}
$$
where $N$ is a fixed positive integer, the $c_i$'s are complex numbers, 
 and only non-negative powers of $z^{1/N}$
occur.  If a finite number of negative powers is allowed, we speak of
fractional Laurent series.  Puiseux's theorem asserts that the set of all
fractional Laurent series is an algebraically closed field. We begin with a
lemma asserting the existence of 
fractional power solutions of a given form for a functional equation.

\begin{lem} \label{Solutions}
Let $h(t)$ be a power series such that $h(0) = 1.$  Then the equation
$$
z = t^k h(t)
$$
has $k$ fractional power series solutions
$\tau_m(z),  \, 0\leq m\leq k-1$ such that 
$$\tau_0(z)= z^{1/k} + \sum_{i=2}^{\infty}  c_i z^{i/k} \qquad \mbox{and}
\qquad  
\tau_m(z)=\xi^mz^{1/k} + \sum_{i=2}^{\infty}  c_i \xi^{mi}z^{i/k}, \ 1\leq
m\leq k-1, $$
where $\xi$ is a primitive $k$-th root of unity. Moreover, if $h(t)$ is
algebraic, then $\tau_0(z),\ldots, \tau_{k-1}(z)$ are also algebraic.
\end{lem}

{\it Proof.}  The proof uses a standard argument.  
%
Suppose that $\tau(z)$ is a fractional
series of the form $\sum_{i \ge 1}  c_i z^{i/k}$; by equating coefficients on
both sides of the equation $z=t^k h(t)$ we show that there are $k$ solutions
of this form.
Expand $\tau(z)^k h(\tau(z))$ as a fractional series.  As $\tau(z)$ has
no constant term,
$$
\tau(z)^k h(\tau(z)) =
c_1^k z + \sum_{j=1}^{\infty} f_j(c_1,c_2, \ldots,c_j,c_{j+1}) z^{1+ j/k},
$$
where $f_j(c_1,c_2, \ldots,c_j,c_{j+1})$ is a linear combination of
terms of the form $c_{r_1}c_{r_2} \ldots c_{r_l}$ with 
$l \geq k$, $1 \leq r_1 \leq r_2 \leq \ldots \leq r_l$, and 
$r_1+r_2+ \cdots +r_l =j+k$. 
The only terms in  $f_j(c_1,c_2, \ldots,c_j,c_{j+1})$ containing
$c_{j+1}$ is $c_1^{k-1} c_{j+1}$, with multiplicity $k$. Hence 
\begin{eqnarray} \label{coeff} 
f_j(c_1,c_2, \ldots,c_j,c_{j+1}) = kc_1^{k-1}c_{j+1} +
\tilde{f}_j(c_1,c_2, \ldots,c_j).  
\end{eqnarray} 
Equating coeffcients on both side of $z=\tau(z)^k h(\tau(z))$,  
we conclude that $c_1$ is a $k$-th root of unity.
Once we have chosen the value of $c_1,$ all the coefficients $c_{j+1}$ can be determined
recursively by the equations $f_j(c_1,c_2, \ldots,c_{j+1}) = 0.$
Thus, there exist $k$ solutions,
$\tau_0(z), \tau_1(z), \ldots, \tau_{k-1}(z),$ where the coefficient
of $z^{1/k}$ in $\tau_m(z)$ equals $\xi^m.$
Suppose that  
$$\tau_0(z)= z^{1/k} + \sum_{i=2}^{\infty}  c_i z^{i/k}.$$
Then the coefficient of $z^{i/k}$ in $\tau_m(z)$ is
$\xi^{mi}c_i$. This follows from the identity 
that
$$
f_j(c_1\xi^m, c_2 \xi^{2m}, \ldots,c_j \xi^{mj}, c_{j+1}\xi^{m(j+1)})
=\xi^{m(j+k)} f_j(c_1,c_2, \ldots,c_j,c_{j+1}) \xi^{m(j+k)}. 
$$


To finish, suppose that $h(t)$ is algebraic. Then so is the series $t^k h(t)$. Let
$Q(t,y)$ be a polynomial such that $Q(t,t^k h(t))=0$. Then $Q(\tau_m(z), z)=0$ and hence, 
$\tau_m(z)$ is algebraic.     
%
\hfill $\Box$

If $h(t)$ is a polynomial, the equation $z=t^k h(t)$ may have fractional
power series solutions with non-zero constant term, but if $h(t)$ has
infinitely many terms, then 
 $\tau_0(z), \tau_1(z), \ldots, \tau_{k-1}(z)$ are
actually all the fractional power series solutions. Of
course the equation may have solutions that are fractional Laurent series. 
Note also that the algebraic equation $P(z,t)$ appearing in the proof 
may have other fractional power
series solutions besides $\tau_0(z),\ldots, \tau_{k-1}(z)$. These solutions
satisfy $z=\tilde{g}(t)$, where $\tilde{g}(t)$ is a solution of $Q(t,y)=0$.
For our
purposes, we will only need existence of the $k$ solutions described
in Lemma~\ref{Solutions}.

\begin{thm} \label{Periodic}
Let $b_0, b_1, \ldots, b_{k-1}$ be a sequence of
non-negative integers of length $k,$ $l = b_{k-1},$
$s_n = ql+b_j,$ where $n = qk+j,$ and
$(\ell_n)$ be a sequence satisfying the Appell relation
$$
\sum_{n=0}^{\infty}  \ell_n t^n \phi(t)^{s_n} = \Psi(t),
$$
where $\phi(t)$ and $\Psi(t)$ are algebraic functions and $\phi(0) = 1.$
Then the ordinary generating function
$ \sum_{n=0}^{\infty} \ell_n z^n $ and the ``section'' generating functions
$\sum_{q=0}^\infty \ell_{qk+j} z^q$, $0\leq j\leq k-1$, 
are algebraic.
\end{thm}

{\it Proof.}
Let $z = t^k \phi(t)^l$ and
$$
Q_j (z) = \sum_{q = 0}^{\infty} \ell_{qk+j} z^q.
$$
Then we can rewrite the Appell relation in the following way:
\begin{eqnarray*}
\Psi(t)  & = &
\sum_{j=0}^{k-1} \sum_{q = 0}^{\infty} \ell_{qk+j} t^{qk+j} \phi(t)^{ql + b_j}
\\
&=& \sum_{j=0}^{k-1}  Q_j(z) t^j \phi(t)^{b_j}.
\end{eqnarray*}

By Lemma \ref{Solutions}, there are $k$ solutions $\tau_i (z)$ of the
form $\xi^i z^{1/k} + \sum_{n = 2}^{\infty} c_n z^{n/k}$, where $\xi$ is a
primitive $k$-th root of unity.  Substituting
the solutions $\tau_i (z)$ into the Appell relation, we get $k$ linear
equations for $Q_j(z):$
$$
\left(\begin{array}{ccccc}
\phi(\tau_0)^{b_0} &  \tau_0\phi(\tau_0)^{b_1} &
\tau_0^2 \phi(\tau_0)^{b_2} & \cdots & \tau_0^{k-1} \phi(\tau_0)^{b_{k-1}}
\\
\phi(\tau_1)^{b_0} &  \tau_1 \phi(\tau_1)^{b_1} &
\tau_1^2 \phi(\tau_1)^{b_2} & \cdots & \tau_1^{k-1}\phi(\tau_1)^{b_{k-1}}
\\
\vdots & \vdots & \vdots & \ddots & \vdots \\
\phi (\tau_{k-1})^{b_0} &  \tau_{k-1} \phi(\tau_{k-1})^{b_1} &
\tau_{k-1}^2  \phi(\tau_{k-1})^{b_2} & \cdots
& \tau_{k-1}^{k-1} \phi(\tau_{k-1})^{b_{k-1}}
\end{array} \right)
\left(
\begin{array}{c}Q_0(z)\\ Q_1(z)\\ \vdots \\Q_{k-1}(z) \end{array}
\right)
=
\left(\begin{array}{c} \Psi(\tau_0) \\ \Psi(\tau_1) \\ \vdots
\\ \Psi(\tau_{k-1}) \end{array}
\right)
$$

The $k \times k$ matrix on the left is non-singular.  Indeed, the
lowest power of $z$ occurring in its determinant is $z^{(k-1)/2}$ and
that power has coefficient equal to the Vandermonde determinant
$\det [\xi^{ij}]_{0 \le i,j \le k-1}.$   Thus, we can solve the
system of linear equations and obtain $Q_j(z)$ as rational functions
in $\tau_0(z), \tau_1(z), \ldots, \tau_{k-1}(z).$  Since the series $\tau_i(z)$ are
algebraic, the functions $Q_j(z)$ and the sum
$ Q_0(z^k) + zQ_1(z^k) + \cdots + z^{k-1}Q_{k-1}(z^k),$ which equals
$\sum_{n=0}^\infty l_n z^n, $ are also algebraic.
\hfill $\Box$

The algebraicity of the generating function for $(a,b)$-paths with an
ultimately periodic boundary follows from the previous theorem.

\begin{cor} \label{PathsAlgebraic}
Let $\underline{s}$ be a ultimately periodic right boundary having height $k$
satisfying the slope condition for $(a,b)$-paths.  Then the ordinary generating function
$$
\sum_{n=0}^{\infty}  \mathrm{SP}_n(a,b, \sigma; \us)  z^n
$$
is algebraic.
In addition, if the initial segment has length $r$, then the section generating functions
$$
\sum_{q=r}^{\infty}  \mathrm{SP}_{qk+j+r}(a,b,\sigma;\us)  z^q,
\qquad j = 0,1,2, \ldots, k-1,
$$ 
are algebraic.
\end{cor}

{\it Proof.} 
Suppose that the boundary $\us$ is the concatenation
$\underline{a}, a_{r-1} + \underline{s'} $, where $\underline{a}$
is the initial segment $a_0,\ldots,a_{r-1},$ and $\underline{s'}$ is a periodic boundary. Then the
Appell relation for $\us$ is
$$\sum_{i=0}^{r-1} \mathrm{SP}_i(a,b,\sigma; \underline{a})t^i
[(1-t)F(t)]^{a_i} 
+ t^r[(1-t)F(t)]^{a_{r-1}}
\sum_{n=0}^{\infty} \mathrm{SP}_{n+r}(a,b,\sigma;\underline{s})t^n [(1-t)F(t)]^{s'_n} =1, 
$$
where $F(t)$ is the series defined in Lemma~\ref{SPAppell}. By Theorem~\ref{Periodic},
 the partial generating function
$\sum_{n=0}^{\infty} \mathrm{SP}_{n+r}(a,b,\sigma;\underline{s})z^n$ as well as
the analogous section
generating functions are algebraic.  Since the
full generating function $\sum_{n=0}^{\infty} \mathrm{SP}_n(a,b,\sigma; \us)$
can be obtained by multiplying the partial generating function by $z^r$ and
adding the polynomial
$\sum_{i=0}^{r-1} \mathrm{SP}_i (a,b,\sigma; \underline{a})z^i,$ the full generating
function is algebraic.
\hfill $\Box$

Corollary \ref{PathsAlgebraic} includes the case of lattice
paths.  Just set $\sigma = 0.$  Thus, we have proved Theorem \ref{TBP}.

As observe earlier, when $\gamma$ is a rational number,
the boundary $\underline{\lceil n/\gamma \rceil + 1}$ is periodic and hence, the
generating functions
$\sum_{n=0}^{\infty} \mathrm{LP}_n(\underline{\lceil n/\gamma \rceil + 1})z^n$
is algebraic.  What can be said about the generating functions
$\sum_{n=0}^{\infty} \mathrm{LP}_n(\underline{\lceil n/\gamma \rceil + 1})z^n$
when $\gamma$ is an algebraic or transcendental number?

\section{The tennis ball problem revisited}

The proof of Theorem \ref{Periodic} provides a recipe for calculating
explicitly the ordinary generating functions.  We show how this might be done
with a concrete example.   Let $\underline{s}$ be the
ultimately periodic boundary
$$\displaystyle 1,
\overbrace{l+1,\ldots,l+1}^{\text{$k$ times}},
\overbrace{2l+1,\ldots,2l+1}^{\text{$k$ times}},
\overbrace{3l+1,\ldots,3l+1}^{\text{$k$ times}},
\ldots\,$$
of height $k$ and width $l.$
The generating function $Q_0(z)=\sum_{q=0}^{\infty} \mathrm{LP}_{qk+1} (\us) z^q $ 
for this boundary and ordinary lattice paths yields a solution
to the generalized tennis ball problem
(posed in \cite{Grim, Merlini} and solved in \cite{deMier}). We note that
there are other choices of boundary, such as
$$
\displaystyle
\overbrace{l+1,\ldots,l+1}^{\text{$k$ times}},
\overbrace{2l+1,\ldots,2l+1}^{\text{$k$ times}},
\overbrace{3l+1,\ldots,3l+1}^{\text{$k$ times}}, \ldots\,$$
that would solve the problem.

Separating the length-$1$ initial segment from the periodic part, 
we have the Appell relation
\begin{eqnarray*}
(1 - t) + \sum_{j=0}^{k-1} Q_j(z) t^{j+1}(1-t)^{l+1} = 1,
\end{eqnarray*}
where $z = t^k (1-t)^l$ and $Q_j(z)=\sum_{q=0}^\infty \mathrm{LP}_{1+qk+j} (\us)
z^q$. Let $\tau_m$, for $0\leq m\leq k-1$, be the $k$ solutions of the equation $z=t^k (1-t)^l$ that
are fractional power series with zero constant term. Using the relation $(1 - \tau_m)^l = z/\tau_m^k,$
we obtain the system of linear equations
$$
\left(\begin{array}{ccccc}
z(1-\tau_0)\tau_0^{1-k} &
z(1-\tau_0)\tau_0^{2-k} &
z(1-\tau_0)\tau_0^{3-k} &
 \cdots &
z(1-\tau_0)
\\
z(1-\tau_1)\tau_1^{1-k} &
z(1-\tau_1)\tau_1^{2-k} &
z(1-\tau_1)\tau_1^{3-k} &
 \cdots &
z(1-\tau_1)
\\
\vdots & \vdots & \vdots & \ddots &\vdots \\
z(1-\tau_{k-1}) \tau_{k-1}^{1-k} &
z(1-\tau_{k-1}) \tau_{k-1}^{2-k} &
z(1-\tau_{k-1}) \tau_{k-1}^{3-k} &
 \cdots &
z(1-\tau_{k-1})
\end{array} \right)
\left(
\begin{array}{c}Q_0(z)\\ Q_1(z)\\ \vdots \\Q_{k-1}(z) \end{array}
\right)
=
\left(\begin{array}{c} \tau_0 \\ \tau_1 \\ \vdots
\\ \tau_{k-1} \end{array}
\right)
$$

We will now solve the system using Cramer's rule.  To do this, we need
to compute two determinants.  First, the determinant of the
$k \times k$ matrix on the left equals
\begin{eqnarray*}
z^k \prod_{j=0}^{k-1}  \tau_j^{1-k}(1-\tau_j)
\det
\left( \begin{array}{ccccc}
1 & \tau_0 & \tau_0^2  & \cdots & \tau_0^{k-1} \\
1 & \tau_1 & \tau_1^2  & \cdots & \tau_1^{k-1} \\
\vdots & \vdots & \vdots & \ddots &
\vdots\\
1 & \tau_{k-1} & \tau_{k-1}^2 & \cdots & \tau_{k-1}^{k-1}
\end{array}\right).
\end{eqnarray*}
Hence, by the formula for the Vandermonde determinant,
this determinant equals
$$
z^k \prod_{j=0}^{k-1}  \tau_j^{1-k}(1-\tau_j)
\, \prod_{0 \leq i< j \leq k-1} (\tau_j-\tau_i).
$$

To compute, say, $Q_0(z)$, also need the determinant
\begin{eqnarray*}
\det \left(\begin{array}{ccccc}
\tau_0 &
z(1-\tau_0)\tau_0^{2-k} &
z(1-\tau_0)\tau_0^{3-k} &
 \cdots &
z(1-\tau_0)
\\
\tau_1 &
z(1-\tau_1)\tau_1^{2-k} &
z(1-\tau_1)\tau_1^{3-k} &
 \cdots &
z(1-\tau_1)
\\
\vdots & \vdots & \vdots & \ddots &\vdots \\
\tau_{k-1}&
z(1-\tau_{k-1}) \tau_{k-1}^{2-k} &
z(1-\tau_{k-1}) \tau_{k-1}^{3-k} &
 \cdots &
z(1-\tau_{k-1})
\end{array} \right).
 \end{eqnarray*}
This determinant equals
 \begin{eqnarray*}
z^{k-1} \prod_{j=0}^{k-1} \tau_j^{2-k}(1 - \tau_j)
 \, \det \left(
\begin{array}{ccccc}
\frac{\tau_0^{k-1}}{1-\tau_0} & 1 & \tau_0 & \cdots & \tau_0^{k-2}
\\
\frac{\tau_1^{k-1}}{1-\tau_1} & 1 & \tau_1 & \cdots & \tau_1^{k-2}\\
\vdots & \vdots & \vdots & \ddots & \vdots\\
\frac{\tau_{k-1}^{k-1}}{1-\tau_{k-1}} & 1 &  \tau_{k-1}  &
\cdots & \tau_{k-1}^{k-2}
\end{array}\right).
\end{eqnarray*}


\begin{lem} \label{Alternant}
 \begin{eqnarray*}
&&\det \left( \begin{array}{ccccc}
\frac{\tau_0^{k-1}}{1-\tau_0} & 1 & \tau_0 & \cdots & \tau_0^{k-2}
\\
\frac{\tau_1^{k-1}}{1-\tau_1} & 1 & \tau_1 & \cdots & \tau_1^{k-2}\\
\vdots & \vdots & \vdots & \ddots & \vdots 
\\
\frac{\tau_{k-1}^{k-1}}{1-\tau_{k-1}} & 1 & \tau_{k-1} &
\cdots & \tau_{k-1}^{k-2}
\end{array}\right)
= (-1)^{k+1}\frac {\prod_{0 \le i < j \le k-1} (\tau_j - \tau_i)}
{\prod_{j=0}^{k-1} (1 - \tau_j) }.
\end{eqnarray*}
\end{lem}
The evaluation of the determinant can be done in many ways, for 
example, as a special case of Theorem (Cauchy$+$) in \cite{Han-Kra} or
Theorem 25 in  \cite{Kratt}, or by the theory of alternants and
symmetric functions, as  described in \cite[Chap. 11]{Muir}.
We include a brief proof here for the convenience of the reader.

{\it Proof.} 
Expanding the  determinant along the first column, we conclude
that it equals the sum
$$
\sum_{j = 0}^{k-1} (-1)^{j+1} \frac{\tau_j^{k-1}}{1-\tau_j}
 \prod_{\genfrac{}{}{0pt}{}{0 \le i < i' \le k-1}{i,i'\neq j}}
 \left({\tau_{i'}}-{\tau_i}\right).
$$
We will obtain this sum in a different way.  From the following
expansion of a singular matrix along the first column,
 \begin{eqnarray*}
0 & = &
\det
\left(\begin{array}{ccccc}
\tau_0^{k-1} & 1 & \tau_0 & \ldots & \tau_0^{k-1}
\\
\tau_1^{k-1} & 1 & \tau_1 & \ldots & \tau_1^{k-1}
\\
\vdots & \vdots & \vdots & \ddots & \vdots
\\
\tau_{k-1}^{k-1}& 1 & \tau_{k-1} & \ldots & \tau_{k-1}^{k-1} \\
1 & 1 & 1 & \ldots & 1 \\
\end{array} \right)
\\
& = &
\left(\sum_{j=0}^{k-1} (-1)^j
\tau_j^{k-1} 
\prod_{\genfrac{}{}{0pt}{}{0 \le i<i'\le k-1}{i,i'\neq j}}
\left( \tau_{i'}-\tau_i \right) \prod_{i\neq j}(1-\tau_i) \right)
+ (-1)^k \prod_{0 \le i < j \le k-1} (\tau_j - \tau_i),
\end{eqnarray*}
we conclude, on dividing both sides by $ \prod_{j=0}^{k-1} (1 - \tau_j), $
that
$$
\sum_{j = 0}^{k-1} (-1)^{j+1} \frac{\tau_j^{k-1}}{1-\tau_j}
 \prod_{\genfrac{}{}{0pt}{}{0 \le i < i' \le k-1}{i,i'\neq j}}
 \left({\tau_{i'}}-{\tau_i}\right)
= (-1)^{k+1} \frac {\prod_{0 \le i < j \le k-1} (\tau_j - \tau_i)}
{\prod_{j=0}^{k-1} (1 - \tau_j) } \,\,.
$$
\hfill $\Box$

Using Cramer's rule and Lemma \ref{Alternant}, we obtain, after some
cancellations,
$$
Q_0(z) =
\frac{(-1)^{k+1}}{z} \prod_{j=0}^{k-1} \frac {\tau_j (z)}{1-\tau_j (z)}.
$$
Thus, we have the following theorem, derived earlier in \cite{deMier}.

\begin{thm} \label{FormulaOne}
Let $\underline{s}$ be the ultimately periodic boundary in the tennis ball
poblem. Then
$$
\sum_{q=0}^{\infty} \mathrm{LP}_{kq+1}(\underline{s}) z^{q}
=
\frac{(-1)^{k+1}}{z} \prod_{j=0}^{k-1} \frac{\tau_j (z)}{1-\tau_j (z)},
$$
where $\tau_0(z),\tau_1(z), \ldots,\tau_{k-1}(z)$ are the $k$ fractional
power series solutions in $z^{1/k}$ with zero constant term
to the equation $z = t^k (1-t)^l.$
\end{thm}

Our formula appears to be different from the one given
in~\cite[Theorem 1]{deMier}.
The formula in that paper is 
$$
Q_0(z)=\frac{-1}{z} \prod_{j=1}^k (1-w_j(z)),
$$
where $w_1(z), w_2(z), \ldots, w_{k}(z)$ are the $k$ fractional
power series solutions of the equation
$(w-1)^k-zw^{k+l}=0.$  We can reorder
the solutions so that $w_{j-1} = 1/(1-\tau_j).$
Thus, the two formulas are equivalent. 

As was done in \cite{deMier}, an explicit formula can be obtained when $k = l.$  
The solution to
the equation $y = t(1-t)$ that does not have constant term is the power series $H(y),$ where
$$
H(y) = \tfrac {1}{2} (1 - \sqrt{1-4y}) = y  - 2y^2 + \ldots \, .
$$
Observing that the equation $z = t^k (1-t)^k$ can be formally rewritten as
to $z^{1/k} = t(1 - t),$ we conclude that
$ H(z^{1/k}), H(\xi z^{1/k}), H(\xi^2 z^{1/k}), \ldots,
H(\xi^{k-1} z^{1/k}),
$
are $k$ fractional power series solutions
to $z = t^k(1-t)^k$ with zero constant term.
Since
$$
\frac {H(y)}{1 - H(y)} =
\frac {1 - \sqrt{1 - 4y}}{1 + \sqrt{1 - 4y}}= \frac{1-\sqrt{1-4y}}{2y}-1,
$$
we obtain the the following theorem.

\begin{cor} \label{FormulaTwo}
Let $\underline{s}$ be the tennis ball boundary with $k = l.$ Then

$$
\sum_{q=0}^{\infty} \mathrm{LP}_{kq}(\underline{s}) z^q
=\frac{(-1)^{k+1}}{z} \prod_{j=0}^{k-1} (C(\xi^j z^{1/k})-1),
$$
where
$$
C(y) = \frac {1 - \sqrt{1 - 4y}}{2y},
$$
the generating function for the Catalan numbers, and
$\xi = e^{2\pi i /k}.$
\end{cor}

When $k = l = 2,$ Corollary \ref{FormulaTwo} specializes (after
simple manipulations) to formulas  
given earlier in \cite{Merlini} and \cite{deMier}.
The simple form of the generating function suggests that after all, there may be
an elementary combinatorial solution to the tennis ball
problem, at least in the case $k=l.$

We can also find formulas for $Q_m(z)$ when $m \ge 1.$   These
generating functions involve more complicated symmetric functions
in $\tau_0(z),\tau_1(z), \ldots, \tau_{k-1}(z).$  For example, when $k = l = 2,$

\begin{eqnarray*}
Q_1(z) & = & \frac{4}{z(1+\sqrt{1+4z})(1+\sqrt{1-4z})} -1
\\
&=&  3 + 22z + 211z^2 + 2306 z^3 + 23270 z^4+ 338444 z^5 + \cdots \,.
\end{eqnarray*}

We remark that the boundaries in the tennis ball problem do not satisfy the
slope condition for $(a,b)$-paths when $a \ge 1.$  There are complicated solutions to the 
analogs of the tennis ball problem for for $(1,b)$-paths.    


\section{Explicit formulas for $(1,b)$-paths}

As remarked earlier, $(1,1)$-paths with the boundary $1,2,3,\ldots \,$ are
in bijection with Schr\"order paths,  which have been studied
in \cite{BSS} and \cite{Mohanty79}.   In this section, we derive
explicit formulas for $(1,b)$-paths with the arithmetic-progression
boundary $c,c+d,c+2d, \ldots\, ,$ where $c$ and $d$ are positive integers.  

We use the Appell relation
$$
\sum_{n=0}^{\infty} \mathrm{SP}_n(1,b,\sigma;\underline{c+id})
\left[t \left( \frac {1-t}{1 + \sigma t^b} \right)^d\right]^n 
=
\left( \frac {1 + \sigma t^b}{1-t} \right)^c
$$
derived in Section 3.  When $b = d = 1,$ then we can solve the quadratic equation
$ z(1+ \sigma t) = t(1-t)$ explicitly.  Choosing the right solution and observing that 
$(1 + \sigma t)/(1-t) = t/z,$  we obtain
\begin{eqnarray*}
\sum_{n=0}^{\infty} \mathrm{SP}_n(1,1,\sigma;\underline{c+i})
z^n
 =  
\left[ \frac {1 - \sigma z - \sqrt{(1 - \sigma z)^2 - 4z^2}}{2z} \right]^c.
\end{eqnarray*}
When $c =1,$  this formula is given
in \cite[p.~44]{BSS}.  We note that the general case (with arbitrary $c$) is
derivable by an easy combinatorial argument from the case $c=1.$
In principle, one can also obtain explicit formulas
when $d= 2$ or $3$ using formulas for solving cubic and quartic equations.

For general $b,$ $c,$ and $d,$ we can obtain explicit formulas for the diagonal-step enumerator
using the Lagrange inversion formula. 
  This formula says that 
$$
n \mathrm{SP}_n(1,b,\sigma;\underline{c+id})
=
[t^{n-1}]
\frac {d}{dt}
\left[\left( \frac {1 + \sigma t^b}{1-t} \right)^c \right]
\left( \frac {1 + \sigma t^b}{1-t} \right)^{nd},
$$
where $[t^{n-1}]$ abbreviates ``the coefficient of $t^{n-1}$ in the series''.
The series on the right side is
$$
c (1 + b \sigma t^{b-1} + (\sigma - b \sigma)t^b)
\frac {(1+ \sigma t^b)^{c+ nd -1}} {(1 -  t)^{c+ nd +1}}.
$$
Using the binomial theorem, we obtain
$$
\mathrm{SP}_n(1,b,\sigma;\underline{c+id})
=
\frac {c}{n}
\left[
T_{n,n-1}(\sigma) +
\sigma  b T_{n,n-b+1}(\sigma) +
\sigma  (1-b) T_{n,n-b}(\sigma)
\right],
$$
where $T_{n,m}(\sigma)$ is the polynomial (with parameters $b,c,d$) defined by
$$
\sum_{m=0}^{\infty}  T_{n,m}(\sigma) t^m = \frac {(1 + \sigma t^b)^{c+nd-1}}{(1-t)^{c+nd+1}}.
$$
Explicitly, 
$$
T_{n,m}(\sigma) = \sum_{j=0}^{\lfloor m/b \rfloor}
\binom {c+nd-1}{j} \binom {c+nd+m-bj}{m-bj} \sigma^j.
$$
In particular, when $b = 1,$ we obtain
$$
\mathrm{SP}_n(1,1,\sigma;\underline{c+id})
=
\frac {c(\sigma + 1)}{n} \sum_{j=0}^{n-1}
\binom {c+nd-1}{j} \binom {c+ (n+1)d -j - 1}{(n-1)-j} \sigma^j.
$$

\section{Parking functions}

Parking functions are rearrangements of lattice paths.
An \emph{ $\us$-parking function of length $n$} is a sequence
$(x_0, x_1, \ldots,x_{n-1})$ of non-negative integers such
that its rearrangement
$(x_{(0)},x_{(1)}, \ldots,x_{(n-1)})$ into
a non-decreasing sequence
satisfies the inequalities: $ x_{(i)} < s_i,$
that is, the non-decreasing rearrangement is an $\us$-lattice path.
In this section, we briefly describe analogs of our results
for parking functions.  We shall freely use results from \cite{KungYan03}.

Let $\mathrm{P}_n(\us)$ be the number of $\us$-parking functions
of length $n.$  Then the following Appell relation holds~\cite[Section 3]{KungYan03}:
$$
1 =\sum_{n=0}^\infty  \mathrm{P}_n(\us) \frac{t^n e^{-s_n t}}{n!}.
$$

\begin{thm} \label{ParkAlgebraic}
Let $\underline{s}$ be a ultimately periodic right boundary.
Then the exponential generating function
$$
\sum_{n=0}^{\infty}  \mathrm{P}_n(\us)  z^n/n!
$$
is algebraic over the extension field $\mathbb{C}(z)\langle T \rangle,$ 
obtained by adding to the field of complex rational functions $\mathbb{C}(z)$   
 the set $T=\{\tau_0(z),\tau_1(z), \ldots,\tau_{k-1}(z)\}$ of the $k$ 
solutions to the equation $z = t^k e^{-lt}$ that are fractional power series
in $z^{1/k}$. 

\end{thm} 

{\it Proof.}  The argument in the proof of
Theorem \ref{Periodic} can be recycled.  The only thing we need
to check that is that $y_m = e^{-\tau_m(z)}$ is algebraic over
$\mathbb{C}(z)\langle T \rangle,$  but this
follows since $y_m^l = z/\tau_m^k.$
\hfill $\Box$

\section{Paths not taken}

The methods in this paper can be applied to other kinds of paths 
(and objects related to paths).  For example, one can consider 
$(a/k,b/k)$-paths with {\it fractional} diagonal steps, that is, 
paths with $E,$ $N,$ and $(a/k,b/k)$ steps, where $a$ and $b$ are 
non-negative integers and $k$ is a positive integer
such that the greatest common divisor of $a,b,k$ is $1.$  When $a = b= 1$
and $k=2,$ then $(a/k,b/k)$-paths to $(n,n)$ with boundary
$1,2,3, \ldots \,$ are in bijection with Motzkin paths from $(0,0)$ to
$(2n,2n)$ as defined in \cite[Ex. 6.38c, p.~238]{EC2}.  The theory
extends naturally to such paths.

It is easy to find boundaries for which the generating function $\sum
\mathrm{LP}_n(\underline{s}) z^n$ is non-algebraic.  By \cite[Chapter
  6]{EC2}, an algebraic power seris is D-finite and the coefficients
of a D-finite power series satisfy a polynomial recursion.  
As $\mathrm{LP}_n(\underline{s})$ is
an increasing sequence of positive integers, there is a polynomial $P(n)$
such that
$$\mathrm{LP}_n(\underline{s}) \le P(n) P(n-1) \cdots P(2)P(1) $$
if $\sum \mathrm{LP}_n(\underline{s}) z^n$ is D-finite.
Thus, if $\underline{s}$ is a
sequence which increases sufficiently rapidly, then the generating function
$\sum \mathrm{LP}_n (\underline{s}) z^n$ is not D-finite, and hence, not
algebraic.  
However, the finer details are elusive.  For example, we do not
have an explicit
example of a boundary with a non-algebraic but D-finite generating
function, although it is reasonable to conjecture that one exists.
It is also unknown whether there are boundaries $\underline{s}$
which are not ultimately periodic with an algebraic generating function
$\sum \mathrm{LP}_n(\underline{s})z^n.$   To conclude with a
concrete problem, is the
generating function for $\mathrm{LP}_n(1,2,2,3,3,3,4,4,4,4, \ldots),$
where there are one $1,$ two $2$'s, three $3$'s, and, in general,
$n$ $n$'s in the boundary, algebraic or D-finite?

\vskip 0.5in

{\large \textbf{Acknowledgement}.}
We thank Marc Noy for many insightful discussions.
Lemma \ref{SPAppell} was found by trial and error with the help of
\textsc{Maple}, the package \texttt{EKHAD} in \cite{a_equals_b}, and
the On-line encyclopedia of integer sequences.  We thank
Christian Krattenthaler for directing us to the paper \cite{Gould} of Gould,
where the lemma was derived earlier.

\makeatletter \renewcommand{\@biblabel}[1]{\hfill#1.}\makeatother

\end{document}